\documentclass[10pt]{article}
\usepackage{epsfig,amssymb,amsmath}
\usepackage{palatino}
\usepackage{verbatim}
\usepackage{enumerate}

\def\baa {\begin{eqnarray*}}
\def\eaa {\end{eqnarray*}}

\def \al {\alpha}
\def \be {\beta}

\def\Frac#1#2{\mbox{\large${\textstyle \frac{#1}{#2}}$}}

\def\Frac#1#2{\mbox{\large${\textstyle \frac{#1}{#2}}$}}

\def \tr {{\rm tr\,}}

\def\A{{\mathbf A}}
\def\B{{\mathbf B}}
\def\E{{\mathbf E}}
\def\PP{{\cal P}}

\def\Ga{\Gamma}

\newtheorem{lemma}{Lemma}[section]
\newtheorem{proposition}[lemma]{Proposition}
\newtheorem{corollary}[lemma]{Corollary}
\newtheorem{theorem}[lemma]{Theorem}
\newtheorem{remark}[lemma]{Remark}

\def\bc  {}
\makeatletter
\@addtoreset{equation}{section}
\makeatother
\def\proof{\medskip\noindent{\bf Proof.} }
\def\qed{\hfill $\Box$}

\newcommand {\ds} {\displaystyle}

\begin{document}

\title{Markov $L_2$-inequality with the Laguerre weight}

\author{G.\,Nikolov, A.\,Shadrin}

\date{}
\maketitle

\begin{abstract}
Let $w_\al(t) := t^{\al}\,e^{-t}$, where $\al > -1$, be the Laguerre
weight function, and let $\|\cdot\|_{w_\al}$ be the associated
$L_2$-norm,
$$
    \|f\|_{w_\al} = \left\{\int_{0}^{\infty} |f(x)|^2
    w_\al(x)\,dx\right\}^{1/2}\,.
$$
By $\PP_n$ we denote the set of algebraic polynomials of degree $\le
n$.
We study the best constant $c_n(\al)$ in the Markov inequality
in this norm
$$
   \|p_n'\|_{w_\al} \le c_n(\al) \|p_n\|_{w_\al}\,,\qquad p_n \in \PP_n\,,
$$
namely the constant
$$
c_n(\al) := \sup_{p_n \in \PP_n}
\frac{\|p_n'\|_{w_\al}}{\|p_n\|_{w_\al}}\,.
$$
We derive explicit lower and upper bounds for the Markov constant
$c_n(\al)$, as well as for the asymptotic Markov constant
$$
c(\al)=\lim_{n\rightarrow\infty}\frac{c_n(\al)}{n}\,.
$$
\end{abstract}

\smallskip


\section{Introduction and statement of the results}
Let $w_\al(t) := t^{\al}\,e^{-t}$, where $\al > -1$, be the Laguerre
weight function, and let $\|\cdot\|_{w_\al}$ be the associated
$L_2$-norm,
$$
    \|f\|_{w_\al} = \left\{\int_{0}^{\infty} |f(x)|^2 w_\al(x)\,dx\right\}^{1/2}\,,
$$
Throughout, $\PP_n$ stands for the set of algebraic polynomials of
degree at most $~n$.
We study here the best constant $c_n(\al)$ in the Markov inequality
in this norm
\begin{equation}\label{e1.1}
\|p_n'\|_{w_\al} \le c_n(\al) \|p_n\|_{w_\al}\,,\qquad p_n \in
\PP_n\,,
\end{equation}
namely the constant
$$
   c_n(\al) := \sup_{p_n \in \PP_n} \frac{\|p_n'\|_{w_\al}}{\|p_n\|_{w_\al}}\,.
$$

Our goal is to obtain \emph{good} and \emph{explicit} lower and
upper bounds for $c_n(\al)$, i.e., to find constants
$\underline{c}(n,\al)$ and $\overline{c}(n,\al)$ such that
$$
\underline{c}(n,\al)\leq c_n(\al)\leq \overline{c}(n,\al)\,,
$$
with a small ratio
$\Frac{\overline{c}(n,\al)}{\underline{c}(n,\al)}$\,. Before
formulating our results here, let us give a brief account on the
results hitherto known.\smallskip

It is only the case $\,\al=0\,$ where the best Markov constant is
known, namely, Tur\'{a}n \cite{pt60} proved that
$$
c_n(0)=\Big(2\sin\frac{\pi}{4n+2}\Big)^{-1}\,.
$$
D\"{o}rfler \cite{pd91} showed that $c_n(\al)=O(n)$ for every fixed
$\,\al>-1\,$ by proving the estimates
\begin{eqnarray}\label{e1.2}
&&c_n(\al)^2\geq\frac{n^2}{(\al+1)(\al+3)}
+\frac{(2\al^2+5\al+6)\,n}{3(\al+1)(\al+2)(\al+3)}
+\frac{\al+6}{3(\al+2)(\al+3)}\,,\label{e1.2}\vspace*{2mm}\\
&&c_n(\al)^2\leq\frac{n(n+1)}{2(\al+1)}\,,\label{e1.3}
\end{eqnarray}
see \cite{pd02} for a more accessible source. In the same paper,
\cite{pd02}, D\"{o}rfler proved for the asymptotic constant
$\,c(\al)=\lim_{\al\rightarrow\infty}\Frac{c_n(\al)}{n}\,$ that
\begin{equation}\label{e1.4}
c(\alpha):=\lim_{n\rightarrow\infty}\frac{c_n(\alpha)}{n}
=\frac{1}{j_{(\alpha-1)/2,1}}\,,
\end{equation}
where $j_{\nu,1}$ is the first positive zero of the Bessel function
$J_{\nu}(z)$\,.

In a recent paper \cite{ns17a} we proved the following
\medskip

\noindent \textbf{Theorem A (\cite[Theorem 1]{ns17a}).} \emph{For
all $\alpha>-1\,$ and $\,n\in \mathbb{N}\,$, $\,n\geq 3\,$, the best
constant $\,c_n(\alpha)\,$ in the Markov inequality
$$
\Vert p^{\prime}\Vert_{w_{\alpha}} \leq c_n(\alpha)\,\Vert
p\Vert_{w_{\alpha}}\,,\qquad p\in\PP_n
$$
admits the estimates
\begin{equation}\label{e1.5}
   \frac{ 2 \big(n+\frac{2\alpha}{3}\big)
         \big(n-\frac{\alpha+1}{6}\big)}
     {(\alpha+1)(\alpha+5)}
< \big[c_n(\alpha)\big]^2 <
  \frac{\big(n+1\big) \big(n+\frac{2(\alpha+1)}{5}\big)}
  {(\alpha+1)\big((\alpha+3)(\alpha+5)\big)^{\frac{1}{3}}}\,,
\end{equation}
where for the left-hand inequality it is additionally assumed that
$n>(\alpha+1)/6$\,.}
\medskip

Clearly, Theorem A implies some inequalities for the asymptotic
Markov constant $\,c(\al)\,$ and, through \eqref{e1.4}, inequalities
for $\,j_{\nu,1}$, the first positive zero of the Bessel function
$\,J_{\nu}\,$ (see \cite[Corollaries~1,\,3]{ns17a}).
\smallskip

We also proved in \cite[Theorem~2]{ns17a} that
$\,c(\al)=O(\al^{-1})\,$, which shows that the upper estimate for
$\,c_n(\al)\,$ in \eqref{e1.5}, though rather good for moderate
$\al$, is not optimal.
\smallskip

Our main result here is an upper bound for $\,c_n(\al)\,$ which is
of the right order with respect to both $n$ and $\al$ as they grow
to infinity.

\begin{theorem}\label{t1.1}
For all $\,n\in \mathbb{N}\,$, $\,n\geq 3\,$, the best constant
$\,c_n(\alpha)\,$ in the Markov inequality \eqref{e1.1} satisfies
the inequality
\begin{equation}
\big[c_n(\al)\big]^2\leq
\frac{4(n+1)\Big(n+3+\frac{3(\al+1)}{4}\Big)}{\al^2+10\al+8}\,,\qquad
\al\geq 2\,.\label{e1.6}
\end{equation}
\end{theorem}
As a consequence of Theorem~\ref{t1.1} and D\"{o}rfler's lower bound
\eqref{e1.2} for $\,c_n(\al)\,$ we show that
$$
c_n^2(\al)\asymp \frac{(n+1)(n+\al+4)}{(\al+1)(\al+8)}\,,\qquad
\al\geq 2\,.
$$
\begin{corollary}\label{c1.2}
For all $\,\al\geq 2\,$ and  $\,n\geq 3\,$ the best constant
$\,c_n(\al)\,$ in the Markov inequality \eqref{e1.1} satisfies
\begin{equation}\label{e1.7}
\frac{(n+1)(n+\al+4)}{2(\al+1)(\al+8)}\leq \big[c_n(\al)\big]^2\leq
\frac{4(n+1)(n+\al+4)}{(\al+1)(\al+8)}\,.
\end{equation}
\end{corollary}
As another consequence, we find the limit value of
$\,(\al+1)c_n^2(\al)\,$ as $\,\al\,$ tends to $-1$, and obtain
asymptotic estimates for $\,\al\,c_n^2(\al)\,$ as $\al\,$ tends to
infinity.
\begin{corollary}\label{c1.3} The best constant
$\,c_n(\al)\,$ in the Markov inequality \eqref{e1.1} satisfies:
\begin{eqnarray*}
(i)\quad \lim_{\al\rightarrow -1}
(\al+1)c_n^2(\al)=\frac{n(n+1)}{2}\,; \vspace*{2mm}\\
(ii)\quad \frac{2n}{3}\leq\lim_{\al\rightarrow \infty}\al\,
c_n^2(\al)\leq 3(n+1)\,.
\end{eqnarray*}
\end{corollary}

Finally, Theorem~\ref{t1.1} provides an upper bound for the
asymptotic Markov constant $\,c(\al)\,$ which is of the correct
order $\,O(\al^{-1})\,$ as $\,\al\,$ tends to infinity. As a
consequence of Theorem~A and Theorem~\ref{t1.1} we have the
following
\begin{corollary}\label{c1.4}
For any $\,\al>-1\,$, the asymptotic Markov constant $\,c(\al) =
\lim\limits_{n\to\infty}n^{-1} c_n(\al)\,$ satisfies the
inequalities
$$
\frac{2}{(\al+1)(\al+5)} < [c(\al)]^2 <
 \begin{cases}\,
\ds{\frac{1}{(\al+1)\sqrt[3]{(\al+3)(\al+5)}}}\,,&\quad -1<\al\leq
 \al^{*}\,,\vspace*{1mm}\\ \, \ds{\frac{4}{\al^2+10\al+8}}\,, &\quad \al>\al^{*}\,,
 \end{cases}
$$
where $\al^{*}\approx 43.4$\,.
\end{corollary}

It is worth noticing here that, for all $\,\al>-1$, the ratio of the
upper and the lower bound for $\,c(\al)\,$ in Corollary~\ref{c1.4}
is less than $\,\sqrt{2}\,$.\smallskip

The rest of the paper is organized as follows. Sect. 2 contains a
brief characterization of the squared best Markov constant
$\,c_n^2(\al)\,$ as the largest eigenvalue of a specific matrix
$\,\A_n\,$. In Sect.~3 we prove some estimates for ratios of Gamma
functions needed for the proof of Theorem~\ref{t1.1}. The proof of
Theorem~\ref{t1.1} is given in Sect. 4. Sect. 5 is concerned with
the evaluation of $\,\|\A_n\|_{F}$, the Frobenius norm of
$\,\A_n\,$, and the bounds for $\,c_n(\al)\,$ implied thereby; in
particular, we reproduce D\"{o}rfler's lower bound \eqref{e1.2}. The
proof of Corollaries~\ref{c1.2}--\ref{c1.4} is given in Sect. 6.
\section{Preliminaries}
It is well-known that the squared best Markov constant $c_n^2(\al)$
equals to the largest eigenvalue of a certain positive definite
matrix $\A_n$. For the reader convenience, here we derive the
explicit form of $\A_n$.
\smallskip

The orthogonal polynomials with respect to the Laguerre weight
function $w_{\al}(x)=x^{\al}e^{-x}$, $x\in \mathbb{R}_{+}$, are
Laguerre polynomials $\{L_m^{(\al)}\}_{m\in \mathbb{N}_0}$, with the
standard normalization
\begin{equation}\label{e2.1}
\Vert
L_m^{\al}\Vert_{w_{\al}}=\Bigg(\frac{\Ga(m+\al+1)}{\Ga(m+1)}\Bigg)^{\frac{1}{2}}
=:\beta_{m+1}\,, \qquad m\in \mathbb{N}_0
\end{equation}
(for the simplicity sake, we suppress the dependance of the $\be$'s
on $\al$).  Further specific properties of the Laguerre polynomials
are (see, e.g., \cite[eqs. (5.1.13), (5.1.14)]{gs75})
\begin{eqnarray}
&&\frac{d}{dx}\{L_m^{(\al)}(x)\}=-L_{m-1}^{(\al+1)}(x)\,,\quad m\in
\mathbb{N}\,,\label{e2.2}\\
&&L_m^{(\al+1)}(x)=\sum_{\nu=0}^{m}L_{\nu}^{(\al)}(x)\,.\label{e2.3}
\end{eqnarray}

Assume that $\hat{p}_n\in\PP_n$, $\,\Vert
\hat{p}_n\Vert_{w_{\al}}=1$, is an extreme polynomial in the $L_2$
Markov inequality \eqref{e1.1}, i.e.,
\begin{equation}\label{e2.4}
\sup\{\Vert p^{\prime}\Vert_{w_\al}^2\,:\,p\in\PP_n\,,\ \Vert
p\Vert_{w_{\al}}=1\} = c_n^2(\al) =
\Vert\hat{p}_n^{\;\prime}\Vert_{w_{\al}}^2\,.
\end{equation}
Without loss of generality, $\,\hat{p}_n$ can be represented in the
form
$$
\hat{p}_n=\sum_{\nu=1}^{n}a_{\nu}\,L_{\nu}^{(\al)}\,,\qquad
a_{\nu}\in \mathbb{R}\,,\quad 1\leq \nu\leq n\,,
$$
then
$$
\Vert
\hat{p}_n\Vert_{w_\al}^2=\sum_{\nu=1}^{n}a_{\nu}^2\beta_{\nu+1}^2=:
\sum_{\nu=1}^{n}t_{\nu}^2=:\Vert \mathbf{t}\Vert^2=1\,,
$$
where $\mathbf{t}=(t_1,\ldots,t_n)^{\top}\in \mathbb{R}^n$ and
$\Vert\cdot\Vert$ is the Euclidean norm in $\mathbb{R}^n$, i.e.,
$\Vert \mathbf{t}\Vert^2=\mathbf{t}^{\top}\mathbf{t}$\,.

By \eqref{e2.1}, \eqref{e2.2} and \eqref{e2.3}, we get
$$
\Vert\hat{p}_n^{\,\prime}\Vert_{w_{\al}}^2\!=\Big\Vert
\sum_{\nu=1}^{n}a_{\nu}\Big(\sum_{\mu=0}^{\nu-1}L_{\mu}^{(\al)}\Big)
\Big\Vert_{w_\al}^2\!\!=\Big\Vert\sum_{\mu=1}^{n}\Big(\sum_{\nu=\mu}^{n}
a_{\nu}\Big)L_{\mu-1}\Big\Vert_{w_\al}^2\!\!=\sum_{\mu=1}^{n}
\Big(\sum_{\nu=\mu}^{n}\frac{\be_{\mu}}{\be_{\nu+1}}\,t_{\nu}\Big)^2\!
=\Vert \mathbf{C}_n \mathbf{t}\Vert^2,
$$
where $\mathbf{C}_n$ is the upper triangular $n\times n$ matrix
$$
\mathbf{C}_n=\begin{pmatrix} \Frac{\be_1}{\be_2} &
\Frac{\be_1}{\be_3} & \cdots & \Frac{\be_1}{\be_{n+1}}\\
0 & \Frac{\be_2}{\be_3} & \cdots & \Frac{\be_2}{\be_{n+1}}\\
\vdots & \vdots & \ddots & \vdots \\
0 & 0 & \cdots & \Frac{\be_{n}}{\be_{n+1}}
\end{pmatrix}\,.
$$
Hence, \eqref{e2.4} admits the equivalent formulation
\begin{equation}\label{e2.5}
c_n^2(\al)=\sup_{\mathop{}^{\mathbf{t}\in \mathbb{R}^n}_{\Vert
\mathbf{t}\Vert=1}}\Vert \mathbf{C}_n \mathbf{t}\Vert^2=
\sup_{\mathop{}^{\mathbf{t}\in \mathbb{R}^n}_{\Vert
\mathbf{t}\Vert=1}} \mathbf{t}^{\top}\mathbf{C}_n^{\top}\mathbf{C}_n
\mathbf{t} =\mu_{\max}(\mathbf{A}_n)\,,
\end{equation}
where $\mu_{\max}(\mathbf{A}_n)$ is the largest eigenvalue of the
positive definite matrix
$\mathbf{A}_n:=\mathbf{C}_n^{\top}\mathbf{C}_n$. A straightforward
calculation reveals that
$$
\mathbf{A}_n=\begin{pmatrix} \frac{\be_1^2}{\be_2^2} &
\frac{\be_1^2}{\be_2\be_3} & \frac{\be_1^2}{\be_2\be_4} & \cdots
& \frac{\be_1^2}{\be_2\be_{n+1}}\\
\frac{\be_1^2}{\be_2\be_3} & \frac{1}{\be_3^2}(\sum_{j=1}^2\be_j^2)
& \frac{1}{\be_3\be_4}(\sum_{j=1}^2\be_j^2) & \cdots &
\frac{1}{\be_2\be_{n+1}}(\sum_{j=1}^2\be_j^2)\\
\frac{\be_1^2}{\be_2\be_4} &
\frac{1}{\be_3\be_4}(\sum_{j=1}^2\be_j^2) &
\frac{1}{\be_4^2}(\sum_{j=1}^3\be_j^2) & \cdots &
\frac{1}{\be_4\be_{n+1}}(\sum_{j=1}^3\be_j^2)\\
\vdots & \vdots & \vdots & \ddots & \vdots\\
\frac{\be_1^2}{\be_2\be_{n+1}} &
\frac{1}{\be_3\be_{n+1}}(\sum_{j=1}^2\be_j^2) &
\frac{1}{\be_4\be_{n+1}}(\sum_{j=1}^3\be_j^2) & \cdots &
\frac{1}{\be_{n+1}^2}(\sum_{j=1}^{n}\be_j^2)
\end{pmatrix}\,.
$$

We observe that the elements $a_{k,i}$ of the matrix $\mathbf{A}_n$
are given by
$$
a_{k,i}=\frac{1}{\be_{i+1}\be_{k+1}}\,\sum_{j=1}^{\min\{k,i\}}\be_j^2
=\begin{cases} \frac{\be_{i+1}}{\be_{k+1}}\Big(\frac{1}
{\be_{i+1}^2}\,\sum_{j=1}^{i}\be_j^2\Big)\,, & i\leq k\,,\vspace*{1mm}\\
\frac{\be_{k+1}}{\be_{i+1}}\Big(\frac{1}
{\be_{k+1}^2}\,\sum_{j=1}^{k}\be_j^2\Big)\,, & i\geq k\,,
\end{cases}
$$
so that
\begin{equation}\label{e2.6}
a_{k,k}=\frac{1}{\be_{k+1}^2}\,\sum_{j=1}^{k}\be_j^2\,,\qquad
a_{k,i}=\begin{cases} \frac{\be_{i+1}}{\be_{k+1}}\,a_{i,i}\,,& i\leq
k\,,\vspace*{1mm}\\ \frac{\be_{k+1}}{\be_{i+1}}\,a_{k,k}\,,& i\geq
k\,.
\end{cases}
\end{equation}
Hence, $\mathbf{A}_n$ can be written in the following simplified
form
\begin{equation}\label{e2.7}
\mathbf{A}_n=\begin{pmatrix} a_{11} & \frac{\be_2}{\be_3}\,a_{11} &
\frac{\be_2}{\be_4}\,a_{11} & \cdots &
\frac{\be_2}{\be_{n+1}}\,a_{11}\\ \frac{\be_2}{\be_3}\,a_{11} &
a_{22} & \frac{\be_3}{\be_4}\,a_{22} & \cdots &
\frac{\be_3}{\be_{n+1}}\,a_{22}\\ \frac{\be_2}{\be_4}\,a_{11} &
\frac{\be_3}{\be_4}\,a_{22} & a_{33} & \cdots &
\frac{\be_4}{\be_{n+1}}\,a_{33}\\ \vdots & \vdots & \vdots & \ddots
& \vdots \\ \frac{\be_2}{\be_{n+1}}\,a_{11} &
\frac{\be_3}{\be_{n+1}}\,a_{22} & \frac{\be_4}{\be_{n+1}}\,a_{33} &
\cdots & a_{nn}
\end{pmatrix}\,.
\end{equation}

We complete this section with giving explicit formulae for $a_{k,k}$
and the trace of $\mathbf{A}_n$.

\begin{proposition}\label{la1}
For every $\,k\in \mathbb{N}\,$ and $\,\al>-1$,
\begin{equation}\label{e2.8}
a_{k,k}=\frac{k}{\al+1}
\end{equation}
and consequently
\begin{equation}\label{e2.9}
\tr (\mathbf{A}_n)=\frac{n(n+1)}{2(\al+1)}\,.
\end{equation}
\end{proposition}

\proof In view of \eqref{e2.6}, we need to show that
\begin{equation}\label{e2.10}
\frac{1}{\be_{k+1}^2}\,\sum_{j=1}^{k}\be_j^2=\frac{k}{\al+1}\,.
\end{equation}
The proof is by induction with respect to $k$. Since
$$
\frac{\be_k^2}{\be_{k+1}^2}=\frac{\frac{\Ga(k+\al)}{\Ga(k)}}
{\frac{\Ga(k+1+\al)}{\Ga(k+1)}}=\frac{k}{k+\al}\,,
$$
\eqref{e2.10} is true for $k=1$. Assuming that \eqref{e2.10} is true
for $\,k-1\in \mathbb{N}$, we obtain
$$
\frac{1}{\be_{k+1}^2}\,\sum_{j=1}^{k}\be_j^2=\frac{\be_k^2}{\be_{k+1}^2}+
\frac{\be_k^2}{\be_{k+1}^2}\,\Big(\frac{1}{\be_{k}^2}\,\sum_{j=1}^{k-1}\be_j^2\Big)
=\frac{k}{k+\al}\,\Big(1+\frac{k-1}{\al+1}\Big)=\frac{k}{\al+1}\,.
$$
Hence, the induction step is performed, and the proof of
\eqref{e2.10} is complete. \qed

\begin{remark}\label{r2.2} D\"{o}rfler's estimate \eqref{e1.3} is
simply the inequality $\;c_n^2(\al)=\mu_{\max}(\A_n)\leq \tr
(\mathbf{A}_n)=\frac{n(n+1)}{2(\al+1)}$\,.
\end{remark}

\section{Estimates for \boldmath{$\frac{\be_i}{\be_k}$}}
We shall need estimates for the elements $\,a_{k,i}$, $\,k\ne i$, of
the matrix $\mathbf{A}_n$ in \eqref{e2.7}, and this requires
estimates for the ratios of the $\be$'s. We prove the following
lemma.
\begin{lemma}\label{l3.1}
For every $\,\al\geq 1\,$ and $\,i,\,k\in \mathbb{N}$, $\,i<k\,$,
there holds
\begin{equation}\label{e3.1}
\frac{\frac{\Ga(i+\al)}{\Ga(i)}}{\frac{\Ga(k+\al)}{\Ga(k)}}\leq
\Bigg(\Frac{i+\Frac{\al-1}{2}}{k+\Frac{\al-1}{2}}\Bigg)^{\al}\,.
\end{equation}
\end{lemma}

\proof It suffices to prove only the case $k=i+1$, for then the
general case will follow from
$$
\frac{\frac{\Ga(i+\al)}{\Ga(i)}}{\frac{\Ga(k+\al)}{\Ga(k)}}=
\prod_{\nu=i}^{k-1}\frac{\frac{\Ga(\nu+\al)}{\Ga(\nu)}}
{\frac{\Ga(\nu+1+\al)}{\Ga(\nu+1)}}\,,
\qquad
\Frac{i+\Frac{\al-1}{2}}{k+\Frac{\al-1}{2}}=\prod_{\nu=i}^{k-1}
\Frac{\nu+\Frac{\al-1}{2}}{\nu+1+\Frac{\al-1}{2}}\,.
$$
Thus, we need to show that
$$
\frac{i}{i+\al}\leq
\Bigg(\Frac{i+\Frac{\al-1}{2}}{i+1+\Frac{\al-1}{2}}\Bigg)^{\al},
\qquad i\geq 1,\ \al\geq 1,
$$
or, equivalently,
\begin{equation}\label{e3.2}
\Bigg(1+\Frac{1}{i+\Frac{\al-1}{2}}\Bigg)^{\al}\leq
1+\frac{\al}{i}\,.
\end{equation}
Clearly, \eqref{e3.2} turns into identity when $\al=1$, so we assume
further that $\al>1$. Set
$$
z=\Frac{1}{i+\Frac{\al-1}{2}}\,,\quad 0<z\leq \frac{2}{\al+1}<1\,,
$$
then
$$
i=\frac{2-(\al-1)z}{2z}\,,
$$
and inequality \eqref{e3.2} becomes
\begin{equation}\label{e3.3}
(1+z)^{\al}\leq 1+\frac{2\al z}{2-(\al-1)z}\,,\qquad 0<z\leq
\frac{2}{\al+1}<1\,,\ \ \al>1\,.
\end{equation}

Assume that $m-1<\al\leq m$, where $m\in \mathbb{N}$, $m\geq 2$. By
Maclaurin's formula, we have
$$
(1+z)^{\al}\leq
1+\sum_{\nu=1}^{m}\frac{\al(\al-1)\ldots(\al-\nu+1)}{\nu!}\,z^{\nu}\,
$$
and it suffices to show that
$$
\sum_{\nu=1}^{m}\frac{\al(\al-1)\ldots(\al-\nu+1)}{\nu!}\,z^{\nu}\leq
\frac{2\al z}{2-(\al-1)z}\,.
$$
Multiplying both sides of this inequality by $2-(\al-1)z>0$ and
arranging the powers of $z$, we arrive at the equivalent inequality
$$
\sum_{\nu=2}^{m+1}\frac{(2-\nu)(\al+1)\al(\al-1)\ldots(\al-\nu+2)}{\nu!}
\,z^{\nu}=:\sum_{\nu=2}^{m+1}a_{\nu}\,z^{\nu}\leq 0\,,
$$
which is obviously true since $z>0$ and $a_{\nu}\leq 0$, $\;2\leq
\nu\leq m+1$.\qed

Lemma~\ref{l3.1} is a particular case of the following more general
statement, which is of independent interest.

\begin{proposition}\label{pa2}
Let $\,i,\,k\in \mathbb{N}$, $\,i<k$.\smallskip

(i)~~~~ If $\,-1<\al\leq 0\,$ or $\,\al\geq 1$, then
\begin{equation}\label{e3.4}
\Big(\frac{i}{k}\Big)^{\al}\leq
\frac{\frac{\Ga(i+\al)}{\Ga(i)}}{\frac{\Ga(k+\al)}{\Ga(k)}}
\leq
\Bigg(\Frac{i+\Frac{\al-1}{2}}{k+\Frac{\al-1}{2}}\Bigg)^{\al}\,.
\end{equation}

(ii)~~~~ If $\,0\leq\al\leq 1$, then
\begin{equation}\label{e3.5}
\Big(\frac{i}{k}\Big)^{\al}\geq
\frac{\frac{\Ga(i+\al)}{\Ga(i)}}{\frac{\Ga(k+\al)}{\Ga(k)}} \geq
\Bigg(\Frac{i+\Frac{\al-1}{2}}{k+\Frac{\al-1}{2}}\Bigg)^{\al}\,.
\end{equation}
\end{proposition}
The proof of Proposition~\ref{pa2} is omitted as we only need its
part given in Lemma~\ref{l3.1}.

\section{Proof of Theorem~\ref{t1.1}}
As was mentioned in Sect. 2, $\,c_n^2(\al)=\mu_{\max}(\A_n)$, where
$\,\mu_{\max}(\A_n)\,$ is the largest eigenvalue of the matrix
$\,\A_n\,$ given by \eqref{e2.7}. It is well-known that
\begin{equation}\label{e4.1}
\mu_{\max}(\A_n)\leq \Vert\A_n\Vert_{*}\,,
\end{equation}
where $\Vert\cdot\Vert_{*}$ is any matrix norm. Here, we shall
exploit $\,\Vert\cdot\Vert_{\infty}$,
$$
\Vert\mathbf{A}_n\Vert_{\infty}=\max_{1\leq k\leq n}
\sum_{i=1}^{n}|a_{k,i}|=\max_{1\leq k\leq n} \sum_{i=1}^{n}a_{k,i}
$$
(notice that $\,a_{k,i}>0$, $\;1\leq i,\,k\leq n$).
Theorem~\ref{t1.1} is an immediate consequence of the following
statement.

\begin{proposition}\label{pa3}
The following inequality holds true:
\begin{equation}\label{e4.2}
\|\A_n\|_{\infty}\leq
\frac{4(n+1)\Big(n+3+\frac{3(\al+1)}{4}\Big)}{\al^2+10\al+8}\,,\qquad
\al\geq 2\,.
\end{equation}
\end{proposition}

We shall need the following lemma, which is proved in \cite{ns17}.

\begin{lemma} \label{l3.3}
Let $\,\al_i>0$, $\gamma_{\min} \le \gamma_i \le \gamma_{\max}$,
$1\le i \le r$, and let
$$
    f(x) := (x+\gamma_1)^{\al_1}(x+\gamma_2)^{\al_2}
              \cdots(x+\gamma_r)^{\al_r}, \qquad
    s := \sum_{i=1}^{r}\al_i\,.
$$
Then, for any $\,x > x_0$,  where $x_0 + \gamma_{\min} \ge 0$, we
have
$$
   \frac{1}{s+1}\, \Big[(t+\gamma_{\min})f(t)\Big]_{x_0}^x
 < \int\limits_{x_0}^{x} f(t)\,dt
 < \frac{1}{s+1}\,(x+\gamma_{\max})f(x)\,.
$$
\end{lemma}

\noindent \textbf{Proof of Proposition~\ref{pa3}.}~ Let us assume
first that $\,\al>2$. For a fixed $k$, $\;1\leq i,\,k\leq n$, we
consider the sum of the elements in the $k$-th row of
$\mathbf{A}_n$,
$$
\sum_{i=1}^{n}a_{k,i}=\sum_{i=1}^{k-1}\frac{\be_{i+1}}{\be_{k+1}}\,a_{i,i}
+a_{k,k}+\sum_{i=k+1}^{n}\frac{\be_{k+1}}{\be_{i+1}}\,a_{k,k}\,.
$$

By Lemma \ref{la1} and Lemma \ref{l3.1}, we have
$$
a_{\nu,\nu}=\frac{\nu}{1+\al}\,,\quad
\frac{\be_{\mu+1}}{\be_{\nu+1}}\leq
\Bigg(\Frac{\mu+\Frac{\al+1}{2}}{\nu+\Frac{\al+1}{2}}\Bigg)^{\frac{\al}{2}}\,,
\quad \mu<\nu\,,\ \ \al\geq 1\,,
$$
hence
\[
\begin{split}
\sum_{i=1}^{n}a_{k,i}&\leq\frac{1}{1+\al}\,\Big[
\Big(k+\frac{\al+1}{2}\Big)^{-\frac{\al}{2}}
\sum_{i=1}^{k-1}i\Big(i+\frac{\al+1}{2}\Big)^{\frac{\al}{2}}+k+
k\Big(k+\frac{\al+1}{2}\Big)^{\frac{\al}{2}} \sum_{i=k+1}^n
\Big(i+\frac{\al+1}{2}\Big)^{-\frac{\al}{2}}\Big]\\
&=:\frac{1}{1+\al}\,\Big[\Big(k+\frac{\al+1}{2}\Big)^{-\frac{\al}{2}}\,S_1
+k+k\Big(k+\frac{\al+1}{2}\Big)^{\frac{\al}{2}}\,S_2\Big]\,.
\end{split}
\]

To obtain an upper bound for $S_1$, we observe that
$\,f_1(x)=x\Big(x+\frac{\al+1}{2}\Big)^{\frac{\al}{2}}\,$ is a
non-negative and increasing function in $(0,\infty)$ to estimate the
sum by an integral, and then apply Lemma~\ref{l3.3} to obtain
$$
S_1\leq\int_{0}^{k}f_1(x)\,dx <
\frac{1}{\frac{\al}{2}+2}\,k\,\Big(k+\frac{\al+1}{2}\Big)^{\frac{\al}{2}+1}
=\frac{2}{\al+4}\,k\,\Big(k+\frac{\al+1}{2}\Big)^{\frac{\al}{2}+1}\,.
$$

Since $\,f_2(x)=\Big(x+\frac{\al+1}{2}\Big)^{-\frac{\al}{2}}\,$ is a
decreasing function in $(0,\infty)$, we estimate $S_2$ from above by
an integral,
$$
S_2\leq\int_{k}^{n+1}f_2(x)\,dx=\frac{2}{\al-2}\,
\Big(k+\frac{\al+1}{2}\Big)^{1-\frac{\al}{2}}\,
\Bigg[1-\Big(\frac{k+\frac{\al+1}{2}}{n+1+\frac{\al+1}{2}}\Big)^{\frac{\al}{2}-1}
\Bigg]\,.
$$

By substituting the above upper bounds for $S_1$ and $S_2$, we
obtain
\begin{equation}\label{e4.4}
\begin{split}
\sum_{i=1}^{n}a_{k,i}&\leq
\frac{k}{\al+1}+\frac{2}{(\al+1)(\al-2)}\,k\Big(k+\frac{\al+1}{2}\Big)
\Bigg[\frac{2(\al+1)}{\al+4}-
\Big(\frac{k+\frac{\al+1}{2}}{n+1+\frac{\al+1}{2}}\Big)^{\frac{\al}{2}-1}
\Bigg]\\
&=:\frac{k}{\al+1}+\frac{2\Big(n+1+\frac{\al+1}{2}\Big)^2}{(\al+1)(\al-2)}
\,\psi_{\al}(k)\varphi_{\al}(y)\,,
\end{split}
\end{equation}
where
$$
\varphi_{\al}(y):=\frac{2(\al+1)}{\al+4}\,y^2-y^{\frac{\al}{2}+1}\,,
\quad y:=\frac{k+\frac{\al+1}{2}}{n+1+\frac{\al+1}{2}}\in (0,1)\,,
$$
$$
\psi_{\al}(k):=\frac{k}{k+\frac{\al+1}{2}}\,,\qquad
$$

For a fixed  $\,\al>2$, the function $\,\varphi_{\al}\,$ has a
unique local extremum in $[0,1]$, a maximum, which is attained at
\begin{equation}\label{e4.5}
y_{\al}=\Big(\frac{8(\al+1)}{(\al+2)(\al+4)}\Big)^{\frac{2}{\al-2}}
=\Big(1-\frac{\al(\al-2)}{(\al+2)(\al+4)}\Big)^{\frac{2}{\al-2}}\in
(0,1)
\end{equation}
and
\begin{equation}\label{e4.6}
\max_{y\in [0,1]}\varphi_{\al}(y)=\varphi_{\al}(y_{\al})
=\frac{2(\al+1)(\al-2)}{(\al+2)(\al+4)}\,y_{\al}^2>0\,.
\end{equation}
We proceed with a further estimation of $y_{\al}^2$. From
\eqref{e4.5} and  $\,\log(1+x)\leq x$, $\,x>-1$, we have
$$
\log
y_{\al}^2=\frac{4}{\al-2}\log\Big(1-\frac{\al(\al-2)}{(\al+2)(\al+4)}\Big)
<-\frac{4\al}{(\al+2)(\al+4)}\,,
$$
hence
$$
y_{\al}^2\leq e^{-\frac{4\al}{(\al+2)(\al+4)}}\leq
\frac{1}{1+\frac{4\al}{(\al+2)(\al+4)}}=\frac{(\al+2)(\al+4)}{\al^2+10\al+8}\,,
$$
where for the last inequality we have used that $\,e^{-x}\leq
\frac{1}{1+x}$, $\;x\geq 0$. Replacing this bound in \eqref{e4.6},
we obtain
$$
\max_{y\in [0,1]}\varphi_{\al}(y)\leq
\frac{2(\al+1)(\al-2)}{\al^2+10\al+8}\,.
$$

This estimate and
$$
\max_{1\leq k\leq n} \psi_{\al}(k)=
\psi_{\al}(n)=\frac{n}{n+\frac{\al+1}{2}}
$$
yields
\[
\begin{split}
\frac{2\Big(n+1+\frac{\al+1}{2}\Big)^2}{(\al+1)(\al-2)} \,\max_{y\in
[0,1]}\,\varphi_{\al}(y)\max_{1\leq k\leq n}\,\psi_{\al}(k)&\leq
\frac{4}{\al^2+10\al+8}\,\frac{n\Big(n+1+\frac{\al+1}{2}\Big)^2}
{n+\frac{\al+1}{2}}\\ &\leq
\frac{4}{\al^2+10\al+8}\,(n+1)\Big(n+1+\frac{\al+1}{2}\Big)\,.
\end{split}
\]

Now we obtain from \eqref{e4.4}
\[
\begin{split}
\sum_{i=1}^{n}a_{k,i}&
<\frac{n+1}{\al+1}+\frac{2\Big(n+1+\frac{\al+1}{2}\Big)^2}{(\al+1)(\al-2)}
\,\max_{y\in [0,1]}\,\varphi_{\al}(y)\max_{1\leq k\leq
n}\,\psi_{\al}(k)\\
&\leq \frac{4}{\al^2+10\al+8}\,(n+1)\Big( n+1+
\frac{\al+1}{2}+\frac{\al^2+10\al+8}{4(\al+1)}\Big)\\
&<\frac{4}{\al^2+10\al+8}\,(n+1)\Big( n+1+
\frac{\al+1}{2}+\frac{\al^2+10\al+9}{4(\al+1)}\Big)\\
&=\frac{4}{\al^2+10\al+8}\,(n+1)\Big( n+3+\frac{3(\al+1)}{4}\Big)\,.
\end{split}
\]
The latter bound is also an upper bound for
$\Vert\A_n\Vert_{\infty}$\,, therefore Proposition~\ref{pa3} is
proved in the case $\,\al>2$.

The proof of the case $\,\al=2\,$ is similar (and somewhat simpler),
and therefore is omitted. \qed

\begin{remark}\label{r4.3} Actually, the above proof works also in
the case $\,1\leq\al< 2\,$ (with a minor modification, e.g.,
$\,\varphi_{\al}\,$ has a minimum instead of maximum in $(0,1)$,
etc.), yielding a similar upper bound for $\,\|\A_n\|_{\infty}$, and
hence for $\,c_n^2(\al)$. However, for small $\,\al\,$ the upper
bound for $\,c_n^2(\al)$ implied by the estimation of
$\,\|\A_n\|_{\infty}$ is worse than the upper bound given in
Theorem~A, and also than the upper bound obtained through the
Frobenius norm of $\,\mathbf{A}_n$.
\end{remark}

\section{The Frobenius norm of \boldmath{$\,\mathbf{A}_n$}}
Let us recall that the Frobenius norm $\,\|\cdot\|_{F}\,$ of a
matrix $\,\B=(b_{i,j})_{n\times n}$ with real elements is defined by
$$
\|\B\|_{F}^2=\sum_{i=1}^{n}\sum_{j=1}^{n}b_{i,j}^2=\tr(\B^{\top}\B)\,.
$$

Since $\,\A_n\,$ is a symmetric and positive definite matrix, we
have
\begin{equation}\label{e5.1}
\|\A_n\|_{F}^2=\tr(\A_n^2)=\mu_1^2+\mu_2^2+\cdots+\mu_n^2\,,
\end{equation}
where $\,0<\mu_1<\mu_2<\cdots<\mu_n=\mu_{\max}(\A_n)\,$ are the
eigenvalues of $\,\A_n\,$, i.e., the zeros of the characteristic
polynomial $\,P_n(\mu)=\det(\mu\E_n-\A_n)$,
$$
P_n(\mu)=\mu^n-b_1\,\mu^{n-1}+b_2\,\mu^{n-2}-b_3\,\mu^{n-3}+\cdots +
(-1)^{n}b_n\,.
$$

As a part of the proof of Theorem~A, in \cite{ns17a} we evaluated
coefficients $\,b_i\,,\ 1\leq i\leq 3$, these coefficients are given
below:
\begin{eqnarray*}
&&b_1=\tr(\A_n)=\frac{n(n+1)}{2(\alpha+1)}\,,\quad
b_2=\frac{(n-1)n(n+1)}{24(\alpha+1)(\alpha+2)(\alpha+3)}
\;\big[3(\alpha+2)n+2(\alpha+6)\big]\,,\vspace*{4mm}\\
&& b_3=\frac{(n-2)(n-1)n(n+1) \big[5(\alpha+2)(\alpha+4)n(n+1)+
8(7\alpha+20)n+12(\alpha+20)\big]}
{240(\alpha+1)(\alpha+2)(\alpha+3)(\alpha+4)(\alpha+5)}\,.
\end{eqnarray*}

Estimates for $\,c_n^2(\al)=\mu_{\max}(\A_n)\,$ are also possible in
terms of solely the first two coefficients, $\,b_1\,$ and $\,b_2\,$.
Indeed, since $\,\tr(\A_n)=b_1\,$ and, by \eqref{e5.1},
$\,\|\A_n\|_{F}^2=b_1^2-2b_2$, we have
$$
b_1-2\,\frac{b_2}{b_1}=\frac{\|\A_n\|_{F}^2}{\tr(\A_n)}\leq\mu_{\max}(\A_n)
\leq \|\A_n\|_{F}=\big(b_1^2-2b_2\big)^{\frac{1}{2}}\,.
$$
Replacing $\,b_1\,$ and $\,b_2\,$ by the expressions above, we
obtain the estimates
\begin{equation}\label{e5.2}
c_n^{4}(\al)\leq b_1^2-2b_2=\frac{n(n+1)}{2(\al+1)^2(\al+3)}\,\Big[
n^2+\frac{2\al^2+5\al+6}{3(\al+2)}\,n+\frac{(\al+1)(\al+6)}{3(\al+2)}
\Big]\,,
\end{equation}
\begin{equation}\label{e5.3}
c_n^{2}(\al)\geq b_1-2\,\frac{b_2}{b_1}= \frac{n^2}{(\al+1)(\al+3)}
+\frac{(2\al^2+5\al+6)\,n}{3(\al+1)(\al+2)(\al+3)}
+\frac{\al+6}{3(\al+2)(\al+3)}\,,
\end{equation}
the second being nothing but the lower estimate \eqref{e1.2} of
D\"{o}rfler.
\smallskip

Slightly weaker but simpler estimates can be obtained on the basis
of \eqref{e5.2} and \eqref{e5.3}.

\begin{proposition}\label{p5.1}
For all $\,n\geq 3\,$,  the best Markov constant $\,c_n(\al)\,$
satisfies the inequalities
\begin{equation}\label{e5.4}
c_n^2(\al)\leq \frac{(n+1)\sqrt{n\big(n+\frac{2(\al+1)}{3}\big)}}
{(\al+1)\sqrt{2(\al+3)}}\,,\qquad \al>-1\,,
\end{equation}
\begin{equation}\label{e5.5}
c_n^2(\al)\geq\begin{cases}
\ds{\frac{n\big(n+\frac{7}{8}\big)}{(\al+1)(\al+3)}}\,,& \al\in
(-1,0)\,, \vspace*{2mm}\\ \ds{\frac{n(n+1)}{(\al+1)(\al+3)}}\,,&
\al\in [0,1]\,,
\vspace*{2mm}\\
\ds{\frac{n\big(n+\frac{2\al+1}{3}\big)}{(\al+1)(\al+3)}}\,, &
\al\geq 1\,.\end{cases}
\end{equation}
\end{proposition}

\proof 1)~~ Inequality \eqref{e5.4} follows from \eqref{e5.2} and
the inequality
$$
n^2+\frac{2\al^2+5\al+6}{3(\al+2)}\,n+\frac{(\al+1)(\al+6)}{3(\al+2)}
\leq (n+1)\Big(n+\frac{2\al+1}{3}\Big)\,.
$$
The latter simplifies to the inequality
$$
\frac{(\al+1)(4n+\al-2)}{3(\al+2)}\geq 0\,,
$$
which is obviously true.\smallskip

2)~~ From \eqref{e5.3} we have
$$
c_n^{2}(\al)\geq
\frac{n\Big(n+\frac{2\al^2+5\al+6}{3(\al+2)}\Big)}{(\al+1)(\al+3)}\
=\frac{n\Big(n+\frac{2\al+1}{3}+\frac{4}{3(\al+2)}\Big)}{(\al+1)(\al+3)}\,,
$$
whence the case $\,\al\geq 1\,$ in \eqref{e5.5} readily follows. The
remaining two cases follow from the observation that
$g(\al)=\frac{2\al^2+5\al+6}{3(\al+2)}$ has a unique local extremum
in $(-1,1]$, a minimum, which is attained at
$\alpha_{*}=\sqrt{2}-2\in (-1,0)$, whence $g(\al)\geq
g(\al_{*})=\frac{4\sqrt{2}}{3}-1>\frac{7}{8}$ for $\al\in (-1,0)$,
and $g(\al)\geq g(0)=1$ for $\al\in [0,1]$. \qed

\begin{remark}\label{r5.1} Estimates \eqref{e5.2} and \eqref{e5.3}
and their consequences \eqref{e5.4} and \eqref{e5.5} are inferior to
the estimates in Theorem~A in the sense that they imply weaker
estimates for the asymptotic Markov constant $\,c(\al)$. In fact, it
can be shown that the upper estimate in Theorem~A is superior to
\eqref{e5.4} for every $\,\al>-1\,$ and $\,n\geq 3$. On the other
hand, for small $\,n\,$ D\"{o}rfler's lower estimate \eqref{e5.3}
and the lower estimates in Proposition~\ref{p5.1} are superior to
the lower estimate in Theorem~A.
\end{remark}
\section{Proof of Corollaries~\ref{c1.2}--\ref{c1.4}}
\noindent \textbf{Proof of Corollary~\ref{c1.2}.} The right-hand
inequality follows from Theorem~\ref{t1.1}: for $\,\al\geq 2\,$ we
have
$$
\big[c_n(\al)\big]^2\leq
\frac{4(n+1)\Big(n+3+\frac{3(\al+1)}{4}\Big)}{\al^2+10\al+8}\leq
\frac{4(n+1)(n+\al+4)}{\al^2+9\al+8}=
\frac{4(n+1)(n+\al+4)}{(\al+1)(\al+8)}\,.
$$
For the left-hand inequality we make use of estimate \eqref{e5.5},
the case $\,\al\geq 1$. For $\,n\geq 3$\, we have
$$
\big[c_n(\al)\big]^2\geq
\frac{n\big(n+\frac{2\al+1}{3}\big)}{(\al+1)(\al+3)}=
\frac{2n\big(n+\al+\frac{n+1}{2}\big)}{3(\al+1)(\al+3)}\geq
\frac{2n(n+\al+2)}{3(\al+1)(\al+3)}>
\frac{2n(n+\al+4)}{3(\al+1)(\al+5)}\,,
$$
where for the last inequality we have used that
$\,f(x)=\frac{x+a}{x+b}\,$ is a decreasing function in $(0,\infty)$
when $\,a>b>0$. A further estimation yields
$$
\frac{2n(n+\al+4)}{3(\al+1)(\al+5)}=\frac{2}{3}\cdot
\frac{\al+8}{\al+5}\cdot \frac{n}{n+1}\cdot
\frac{(n+1)(n+\al+4)}{(\al+1)(\al+8)}> \frac{2}{3}\cdot\frac{3}{4}
\cdot\frac{(n+1)(n+\al+4)}{(\al+1)(\al+8)},
$$
which proves the left-hand inequality in Corollary~\ref{c1.2}.\qed
\medskip

\noindent \textbf{Proof of Corollary~\ref{c1.3}.}~ (i) From
\eqref{e5.3} we deduce
$$
\lim_{\al\rightarrow -1} (\al+1)c_n^2(\al)\geq \frac{n(n+1)}{2}\,,
$$
while from the upper estimate in Theorem~A we obtain
$$
\lim_{\al\rightarrow -1} (\al+1)c_n^2(\al)\leq \frac{n(n+1)}{2}
$$
(notice that the same conclusion follows from \eqref{e5.4}).
\smallskip

\noindent (ii) The right-hand inequality follows from
Theorem~\ref{t1.1}, and the left-hand inequality follows from
\eqref{e5.5}. \qed\medskip

\noindent \textbf{Proof of Corollary~\ref{c1.4}.}~ The lower
estimate is a consequence from Theorem~A, while the upper estimates
follow from Theorem~A and Theorem~\ref{t1.1}, respectively\,.\qed

\bigskip\bigskip

\noindent {\bf Acknowledgement.} The research of the first-named
author is supported by the Bulgarian National Research Fund through
Contract DN 02/14.


\newpage

\noindent
{\sc Geno Nikolov} \smallskip\\
Department of Mathematics and Informatics\\
Universlty of Sofia \\
5 James Bourchier Blvd. \\
1164 Sofia \\
BULGARIA \\
{\it E-mail:} {\tt geno@fmi.uni-sofia.bg}

\bigskip\bigskip\noindent
{\sc Alexei Shadrin} \smallskip\\
Department of Applied Mathematics and Theoretical Physics (DAMTP) \\
Cambridge University \\
Wilberforce Road \\
Cambridge CB3 0WA \\
UNITED KINGDOM \\
{\it E-mail:} {\tt a.shadrin@damtp.cam.ac.uk}

\end{document}